\documentclass[12pt]{amsart}
\usepackage{a4,amsmath,amssymb,amsthm,eucal}
\usepackage[english]{babel}
\usepackage{graphicx}
\usepackage{xcolor}
\usepackage{epsfig}

\usepackage{color}

\usepackage[normalem]{ulem}

\definecolor{blue-violet}{rgb}{0.54,0.17,0.89}
\definecolor{amethyst}{rgb}{0.6,0.4,0.8}
\definecolor{darkviolet}{rgb}{0.58, 0.0, 0.83}
\definecolor{darkgreen}{rgb}{0,.4,0}
\definecolor{mixedgreen}{rgb}{0.3,0.6,00}
\definecolor{bananayellow}{rgb}{1.0, 0.88, 0.21}
\definecolor{arylideyellow}{rgb}{0.91, 0.84, 0.42}
\definecolor{bananamania}{rgb}{0.98, 0.91, 0.71}

\newcommand{\REMOVEsk}[1]%
           {{\color{magenta}\sout{#1}}}

\theoremstyle{remark}

\numberwithin{equation}{section}

\begin{document}
\title{Uniform solvability for families of linear systems on time scales}
\author{Sergey Kryzhevich}
\address[Sergey Kryzhevich]{
Institute of Applied Mathematics, Faculty of Applied Physics and Mathematics, Gda\'nsk University of Technology, 80-233 Gdańsk, Poland
\and
BioTechMed Center, Gda\'nsk University of Technology, 80-233 Gda\'nsk, Poland }
\email[Sergey Kryzhevich]{serkryzh@pg.edu.pl}

\thanks{The work of the first author was supported by Gda\'{n}sk University of Technology by the DEC 14/2021/IDUB/I.1 grant under the Nobelium - 'Excellence Initiative - Research University' program.}

\date{March 17, 2022}

\begin{abstract}
	We give explicit criteria of solvability for families of linear systems on time scales. We introduce a new method of embedding a time scale into a non-autonomous system of ODEs. This will be the first step to implementing the structural stability result obtained by one of the co-authors together with V.\, A.\, Pliss to time scale dynamics.
\end{abstract}


\maketitle

\bigskip

\section{Introduction}

Non-autonomous systems of ordinary differential equations are usually regarded as more sophisticated compared to classical flows engendered by autonomous systems of ordinary differential equations. The principal trouble is that periodic (and even recurrent) solutions do not play an important role in the dynamics of non-autonomous systems unless the right-hand side of the studied system is periodic. 

For flows, engendered by autonomous systems of ODEs, the structural stability, and $\Omega$ -- stability theorems were obtained by Robinson \cite{R75} and Ma\~n\'e \cite{M79}. Both of them are related to the so-called Axiom A which is also the core assumption for various further results in Stability Theory. 

For Axiom-A flows, the non-wandering set can be represented as a disjoint union of a finite number of closed transitive invariant subsets. For any of these subsets, there is a neighborhood where the behavior of solutions is hyperbolic (dimensions of stable and unstable manifolds can, however, vary). 

Moreover, all solutions of Axiom-A systems spend a finite (and uniformly bounded) amount of time out of the union of those neighborhoods. In other words, one can say that linear approximations are hyperbolic on big intervals of time, depending on a particular solution. 

The above conditions can be reformulated for non-autonomous systems, see \cite{KP03} and \cite{P77}-\cite{P81} where sufficient conditions for structural stability were obtained. However, it is still unclear if the mentioned conditions are close to necessary ones.  

We are going to translate these results to the language of time scale dynamics. This is a relatively new area in dynamical systems, first introduced by \cite{AH90}. Firstly, this can be regarded as a combination of discrete and continuous dynamics and, secondly, as a generalization of numerical methods with non-uniform steps. A survey on the theory of such systems including Stability Theory is given in \cite{BP03} and the relatively recent book \cite{M16}, see also references therein. However, no global structural stability conditions for generic time scale systems have been established yet. It should be highlighted that we do not just copy-paste the result from the ODE theory. Some additional conditions respecting the geometry of the time scale have to be added. For instance, we consider a particular case of a non-syndetic time scale (i.e. that with an unbounded graininess function).  

The main objective of this paper is to prove analogs of mentioned results for time scale systems. For this purpose, we introduce a renormalization of time scales and, re-define the concept of hyperbolicity.  for time scale systems. Moreover, we show how linear systems on time scales can be reduced to ordinary differential equations.  

For the sake of convenience, we deal with a dynamical system on a Euclidean space (unlike the classical case of ordinary differential equations, in the case of non-periodic time scales there is no principle difference between autonomous and non-autonomous systems). 

The paper is organized as follows.

In Section 2 we recall the results of the classical structural stability theory. Next, in Section 3, we recall the basic concepts of time scale dynamics. Section 4 is the core of the paper. There, we discuss how a linear system on a time scale can be represented as a system of ordinary differential equations (and also, how the concept of hyperbolicity can be translated into the language of time scale dynamics). Finally, in Section 5, we formulate the structural stability conjecture which will be the main objective of further research.

In a nutshell, we follow the main idea of the quoted papers by V.\, A.\, Pliss. First of all, we obtain some conditions on uniform solvability for families of linear systems. Thereafter, we are going to develop Perron's approximations for time scale systems to be able to find solutions for nonlinear systems. 

\section{Structural stability for systems of ODEs.}

Let $V$ be a $C^1$ -- smooth vector field on a $C^1$ -- smooth Riemannian manifold $M$, and $\Omega_V$ be the non-wandering set of the considered autonomous system of ordinary differential equations
$$\dot x=V(x).$$

\bigskip

We recall the famous Axiom $A'$ for flows/vector fields that claims that 
\begin{enumerate}
\item $\Omega_V$ can be represented as a union of two disjoint closed sets $\Omega^1_V$ and $\Omega^2_V$ where $\Omega^1_V$ consists of finitely many isolated fixed points and $\Omega^2_V$ 
does not contain any stationary point;
\item the set $\Omega_V$ is hyperbolic.
\item periodic orbits are dense in $\Omega^2_V$.
\end{enumerate}

Then the following statement is well-known in Hyperbolic Theory.

\bigskip

\noindent \textbf{Smale's spectral decomposition theorem}.
Let the vector field $V$ satisfy Axiom $A'$. Then the non-wandering set $\Omega_V$ admits the unique representation:
$$\Omega_f=\Omega_1\bigcup \Omega_2\bigcup \ldots \bigcup \Omega_N$$
with all components being closed, disjoint, invariant, and transitive.

\bigskip

Structural stability requires an additional condition:

\bigskip

\noindent\textbf{Geometric strong transversality condition}.
An $A'$ flow, defined by a vector field $V$, satisfies the strong transversality condition if for any $p,q \in \Omega_V$ the manifolds $W^s(p)$ and $W^u(q)$ intersect transversally.

\bigskip

Now let us recall the famous result by C.\, Robinson \cite{R75}.

\bigskip

\noindent\textbf{Structural Stability Theorem.} \emph{For $C^1$ -- smooth vector fields, Axiom $A'$ and Geometric Strong Transversality Condition imply structural stability.}

\bigskip

In this section we recall a 'non-autonomous' analog of the above result. Now let us consider a non-autonomous system
\begin{equation}\label{e1.1}
\dot x=V(t,x), \qquad t\in {\mathbb R}, \quad x\in {\mathbb R}^n; \end{equation}
where
$$|V(t,x)|\le M, \qquad |D_xV(t,x)|\le M.$$
Let 
$x(t,t_0,x_0)$ be a solution of Eq.\, \eqref{e1.1} corresponding to initial conditions $x(t_0)=x_0$.

Let us also consider a perturbed system:
\begin{equation}\label{e1.2}
\dot x=V(t,x)+Y(t,x), \qquad |Y(t,x)|\le \delta, \qquad |D_xY(t,x)|\le \delta. 
\end{equation}
with $y(t,t_0,x_0)$ being a solution corresponding to the same initial conditions.

Let us recall the main definitions of the Hyperbolic Theory for non-autonomous systems of ordinary differential equations.

Consider a linear system
\begin{equation}\label{eqlin}
\dot x=A(t)x, \qquad t\in [t_-,t_+], 
\end{equation}
Let $\Phi(t,t_0)=\Phi(t)\Phi^{-1}(t_0)$ be the Cauchy matrix of that system.

We say that system \eqref{eqlin} is hyperbolic on the segment $[t_-,t_+]$ if for any $t\in [t_-,t_+]$ there exists a decomposition
$${\mathbb R}^n=E^s(t)\oplus E^u(t): \qquad \Phi(t,\tau)E^{s,u}(\tau)=E^{s,u}(t)$$
such that
$$\begin{array}{c}
|\Phi(t,\tau) x_0|\le a\exp(-\lambda(t-\tau))|x_0|, \qquad  
t\ge \tau, x_0\in E^s(\tau);\\
|\Phi(t,\tau) x_0|\le a\exp(\lambda(t-\tau))|x_0|, \qquad
t\le \tau, x_0\in E^u(\tau).
\end{array}$$
for some $a>0$ and $\lambda>0.$

Now we consider the set of linear systems
\begin{equation}\label{e1.3}
\dot x = A(t,x_0)x, \qquad x_0 \in {\mathbb R}^n, \qquad
A(t,x_0) = \dfrac{\partial V(t,x(t,0,x_0))}{\partial x} 
\end{equation}

We make the following assumptions.
Let $a_0 > 0$, $\lambda > 0$, $\alpha > 0$, be such that for any $x_0\in {\mathbb R}^n$:
\begin{enumerate}
\item There exists a number $a\in (0,a_0)$ and instants of time 
$$-\infty= \tau_0 < \tau_1< \ldots < \tau_s < \tau_{s+1} = \infty, \qquad 0 < s < n$$ 
such that on any of the segments $[\tau_j,\tau_{j+1}]$, $j = 0,\dots, s$ system \eqref{e1.3} is hyperbolic with constants $a$ and $\lambda$. Let $E^s_j$ and $E^u_j$ be the corresponding stable and unstable spaces. 
\item $\dim E_j^s(\tau_{j+1}) < \dim E_{j+1}^s(\tau_{j+1}), \qquad j = 0, \ldots s-1.$
\item Spaces $E_j^u(\tau_{j+1})$ and $E_{j+1}^s(\tau_{j+1})$ intersect transversally and angles between them satisfy inequalities 
$\sphericalangle E_j^u(\tau_{j+1}), E_{j+1}^s(\tau_{j+1})\ge \alpha$.
\end{enumerate}

It was proved \cite{KP03} that such 'uniform hyperbolicity' of family Eq.\, \eqref{e1.3} with appropriate implies structural stability for nonlinear system \eqref{e1.1}.  

Let the function $T(a, \lambda, \alpha)$ satisfy inequalities
$$
36a^2\exp\left(\dfrac{-\lambda T}{3}\right)<\dfrac{\alpha}{8}\sin \dfrac{\alpha}{4}, \qquad
3a\left(\dfrac{2}{\sin(\alpha/2)}+1\right)\exp\left(\dfrac{-\lambda T}{3}\right)<1.
$$

\bigskip

\noindent\textbf{Theorem 1 \cite{KP03}.}\emph{
Let the system \eqref{e1.3} satisfy the above conditions I-III and for any $j$ the conditions
$$\tau_{j+1}-\tau_j> T(a,\lambda,\alpha), \qquad j = 0,\ldots s.$$
Then for any $\varepsilon>0$ there exists a $\delta>0$ such that if condition \eqref{e1.2} is satisfied, there exists a homeomorphism $\varphi$ of the space ${\mathbb R}^n$ such that
$$|x(t,0,x_0) - y(t,0,\varphi(x_0))| <\varepsilon$$
for any $x_0 \in {\mathbb R}^n$.}

\bigskip

Our main objective (not limited to that paper) is to translate the mentioned result into the language of time scale dynamics.

\section{Timescale dynamics}

In this section, we recall some basic definitions of time scale dynamics.

\bigskip

\noindent\textbf{Definition 1.} 
We say that a subset ${\mathbb T}\subset {\mathbb R}$ is a \emph{timescale} if it is closed and 
$$\inf {\mathbb T}=-\infty, \qquad \sup {\mathbb T}=+\infty$$
We always assume for convenience that $0\in {\mathbb T}$.

\bigskip

Given a time scale, we consider the  \emph{forward jump operator}
$$\sigma: {\mathbb T} \to {\mathbb T}: \sigma(t)=\inf\{s\in {\mathbb T}: s>t\}$$
and consider the so-called \emph{graininess function}, see Fig.~1:
$$\mu(t):=\sigma(t)-t.$$  
We say that a point $t\in {\mathbb T}$ is \emph{right-dense} if 
$\mu(t)=0$ and \emph{left-dense} is $\sigma (s)=t$ implies $s=t$.
We say that a time-scale is \emph{syndetic} if 
$\sup_{t\in {\mathbb T}} \mu(t)<+\infty$.

\begin{figure}[!ht]
\begin{center}
\includegraphics[height=0.8in]{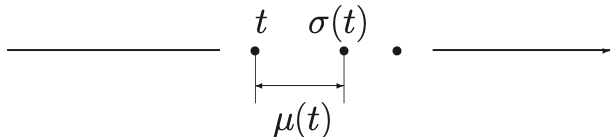}
\caption{\footnotesize A sample of a time scale.}
\end{center}
\end{figure}

The classic concepts of Analysis may also be defined for time scales.

\bigskip

\noindent\textbf{Definition 2.}
A function $f:{\mathbb T} \mapsto {\mathbb R}$ is called \emph{rd-continuous} provided it is continuous at right-dense points in ${\mathbb T}$ and finite left-sided limits exist at left-dense points in ${\mathbb T}$. Denote the class of rd-continuous functions by ${\mathcal C}_{rd} = {\mathcal C}_{rd}({\mathbb T},{\mathbb R})$.

\bigskip

\noindent\textbf{Definition 3.}
The function $f:{\mathbb T} \mapsto {\mathbb R}$ is called {\it $\Delta$-differentiable} at a point $t\in {\mathbb T}$ if there exists $\gamma \in {\mathbb R}$ such that for any $\varepsilon > 0$ there exists a neighborhood $W$ of $t$ satisfying
$$|[f(\sigma(t)) - f(s)] - \gamma[\sigma(t) - s]| \le \varepsilon|\sigma(t) - s|$$
for all $s\in W$. In this case we write $f^\Delta(t)=\gamma.$

\bigskip

For example, when ${\mathbb T} = {\mathbb R}$, $f^\Delta(t) =\dot f(t)$. When ${\mathbb T} = {\mathbb Z}$, $f^\Delta(t)$ is the standard forward difference operator $f(n + 1) - f(n)$.

In time scale dynamics, there is a dual concept of the so-called $\nabla$ - derivative, corresponding to the {\it left jump} operator. We do not consider such a derivative in this paper. Probably, this will be done in the future. 

Also, we define integrals over time scales.

\bigskip

\noindent\textbf{Definition 4.} If $F^\Delta(t) = f(t)$, $t\in {\mathbb T}$, then $F$ is a \emph{$\Delta$ - antiderivative} of $f$, and the Cauchy {$\Delta$ - integral} is given by formula 
$$\int_{a}^b f (t)\, \Delta t = F (b) - F (a),\quad \mbox{for all}\quad a, b \in {\mathbb T}.$$

\bigskip

We recall a simple formula which relates integrals over time-scale with classic ones (this formula can also be used as a definition of integrals and derivatives over time scales)
$$\int_{a}^b f (t)\, \Delta t=\int\limits_{[a,b]\cap {\mathbb T}} f(t) d\mathbf{m}(t)+\sum\limits_{\mu(t_i)>0} f(t_i)\mu(t_i).$$

\bigskip

\noindent\textbf{Definition 5.} A function $p : {\mathbb T} \mapsto {\mathbb R}$ is called \emph{regressive} provided that $$1 + \mu(t)p(t)\neq 0$$ for all $t\in {\mathbb T}$ and \emph{positively regressive} if $1+\mu(t)p(t)>0$ for all $t \in {\mathbb T}$. A matrix mapping ${\mathcal A} : {\mathbb T} \mapsto M^n({\mathbb R})$ is called \emph{regressive} if for each $t\in {\mathbb T}$ an $n\times n$ matrix $E_n + \mu(t){\mathcal A}(t)$ is invertible, and \emph{uniformly regressive} if in addition the matrix function $(E_n + \mu(t){\mathcal A}(t))^{-1}$ is bounded.

\bigskip

Regressivity is equivalent to the backward uniqueness of solutions of the system $x^\Delta=A(t) x$. 

We give a counterexample of a scalar equation, which is not regressive: $$x^\Delta=-x,\quad {\mathbb T}={\mathbb N}.$$ This equation has only one solution, defined on the whole time scale ${\mathbb T}$ that is $x\equiv 0$.

For a linear system
\begin{equation}\label{e1}
x^\Delta=A(t)x, \qquad t \in {\mathbb T},
\end{equation}
stability and asymptotic stability are defined similarly to ODEs. 

The following result is an analog of that of the theory of linear ODEs.

\bigskip

\noindent\textbf{Theorem (Choi and DaCunha)}, \cite{CIK08},  \cite{D05}. 
Let an $rd$ - continuous matrix function $A$ be regressive. Linear system \eqref{e1} is stable if and only if all its solutions are bounded on ${\mathbb T}\bigcap [0,+\infty)$. It is uniformly stable if and only if there exists a positive constant $\gamma$, such that $|\Phi_A(t,t_0)|\le \gamma$, $\forall 
t_0,t \in {\mathbb T}, t\ge t_0$.

\bigskip

\noindent\textbf{Definition 6.} For $p \in {\mathcal R}$, we define the generalized exponential function $e_p(t, s)$ by
$$ e_p(t,s) = \exp\left(\int_s^t \xi_{\mu(\tau)}p(\tau)\, \Delta \tau\right).$$
where $\xi_h$ is the cylinder transformation given by formula $\xi_h(z)=\log(1+zh)/h$ if $h\neq 0$, $\xi_h(z)=z$ if $h=0$.

This is the solution of equation $x^\Delta=p(t)x$ with initial conditions $x(s)=1$.

\bigskip

Other techniques, e.g. Gr\"{o}nwall–Bellmann Inequality, Lyapunov exponents theory can be developed for systems on time scales.

\bigskip

However, stability, Lyapunov exponents, and even uniqueness of solutions depend significantly on a time scale. This is a serious obstacle to applying the classic Hyperbolic Theory for time scale dynamics. For example, if we consider ${\mathbb T}=\{2^n: n\in {\mathbb N}\}$, for any bounded matrix $A(t)$ all solutions have non-positive Lyapunov exponents.

\section{Hyperbolicity for time scale systems}

\noindent\textbf{1. Renormalization}

Let ${\mathbf m}$ be the Lebesgue measure on the real line. Recall that $\mu(t)$ is the graininess function of the time scale. We fix 
$$t_0=\inf \{{\mathbb T}\bigcap [0,\infty)\}$$ 
and define the function 
$s:{\mathbb R} \to {\mathbb R}$ as follows: 
\begin{enumerate}
\item $s(t)={\mathbf m}([t_0,t])+\sum_{\tau\in [t_0,t):\mu(\tau)>0} \ln(1+\mu(\tau))$ if $t\in {\mathbb T}, t\ge t_0$.
\item $s(t)=s(t_1)+\ln (1+t-t_1)$ if $t\notin {\mathbb T}, t\ge t_0$. Here $t_1=\max({\mathbb T}\bigcap (-\infty,t])$.
\item $s(t)=-{\mathbf m}([t_0,t])-\sum_{\tau\in [t_0,t):\mu(\tau)>0} \ln(1+\mu(\tau))$ if $t\in {\mathbb T}, t< t_0$.
\item $s(t)=s(t_1)-\ln (1+t_2-t)$ if $t\notin {\mathbb T}, t\ge t_0$. Here $t_2=\min({\mathbb T}\bigcap [t,\infty))$.
\end{enumerate}

Evidently $s(t_0)=0$, $s(t)$ is strictly monotonous, continuous, unbounded in both directions and $|s(t)|\le |t-t_0|$ for any $t \in {\mathbb R}$. Given a time scale $\mathbb T$ we introduce the renormalized time scale 
${\mathbb S}:=s(\mathbb T)$, Fig.\, 2. 

\begin{figure}[!ht]
\begin{center}
\includegraphics[height=3in]{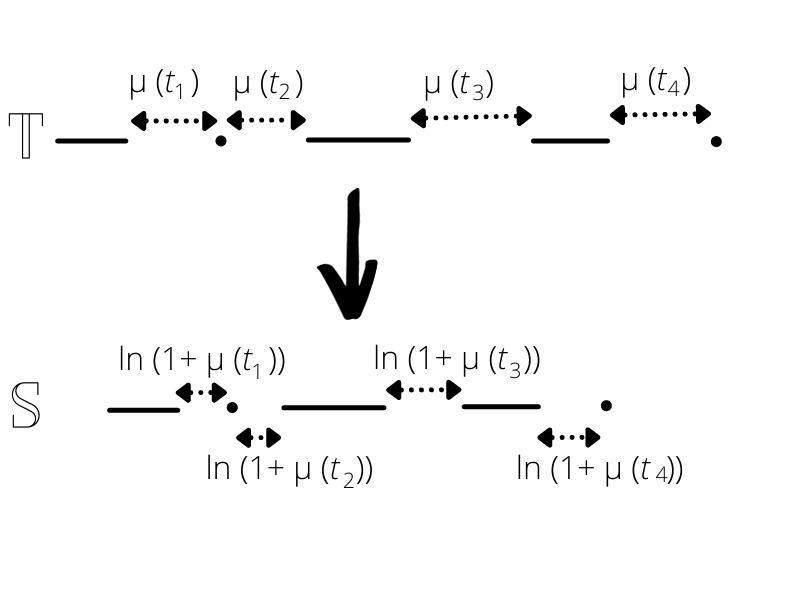}
\caption{\footnotesize A sample of renormalization procedure.}
\end{center}
\end{figure}

\bigskip

\noindent\textbf{2. Equivalent linear system}.
For any system
\begin{equation}\label{e7}
x^\Delta =A(t) x +f(t)
\end{equation}
on the time scale $\mathbb T$ with $A$ and $f$ being bounded uniformly $rd$-continuous there exists a system 
\begin{equation}\label{e8}
\dot x =B(t) x +g(t) 
\end{equation}
such that $\Phi_A(t,\tau)=\Phi_B(s(t),s(\tau))$ for all $t,s\in {\mathbb T}$. Here $\Phi_A$ and $\Phi_B$ are fundamental matrices of systems \eqref{e7} and \eqref{e8}, respectively.  Besides, given the right-hand side $f$ we construct the function $g$ so that for any solution $\varphi (t)$ of system \eqref{e7} there exists a solution $\psi(t)$ of Eq.\, \eqref{e8} such that $\varphi(t)=\psi(s(t))$ for any $t\in {\mathbb T}$.

Here we may assume that $B$ is constant on all connected components of the completion of ${\mathbb T}$.

\bigskip

\noindent\textbf{Definition 7.} Let $(t_1,t_2)$ be an interval such that $t_{1,2} \in {\mathbb T}$, $(t_1,t_2)\bigcap {\mathbb T}=\emptyset$.

\bigskip

\noindent\textbf{Lemma 1.} \emph{The matrix $B$ and the vector function $f$ can be found by the following rule.
\begin{enumerate}
\item $A(t)=B(s(t))$, $f(t)=g(s(t))$ for all $t\in {\mathbb T}$.
\item Let $t\notin {\mathbb T}$. We consider $t_0=\sup_{\tau\in {\mathbb T}: \tau\le t}$  and define 
\begin{equation}\label{e9}
\begin{array}{c}
B(s(t))= \dfrac{{\mathrm{Ln}}\, (E+\mu(t_0)A(t_0))}{\ln(1+\mu(t_0))},\\[5pt] g(s(t))=\dfrac{A^{-1}(t_0){\mathrm{Ln}}\,(E+\mu(t_0)A(t_0))}{1+\mu(t_0)}f(t)
\end{array}
\end{equation}
\end{enumerate}}

\bigskip

\noindent\textbf{Remark 1.} Here $\mathrm{Ln}$ stands for the matrix logarithm. This function is multi-valued and it could be imaginary even if the matrix $A$ is real. However, it can be canonically defined as a real value provided 
$\sup_t\mu(t)\|A(t)\|<1$ where the matrix norm is taken in an appropriate basis. 
For this we can use the classic Taylor decomposition
$$\mathrm{Ln}\,(E+\mu A)=\sum_{k=1}^\infty (-1)^{k+1} \mu^k A^k.$$

Meanwhile, the second formula of \eqref{e9} can be correctly defined even in the case when $A(t_0)$ is non-invertible. This follows from the fact that the function $\ln(1+\mu x)/x$ is holomorphic in a neighborhood of zero.

\bigskip

The latter expression, being well-defined for small values of $\mu$ and a fixed matrix $A$ can be extended for all positive values of $\mu$ unless the matrix $A$ has negative eigenvalues.

If matrix $A$ has negative eigenvalues we still link a negative value $-\mu$ with zero going around the 'problematic' values $\lambda_j$ if any from above. Here $\lambda_j$ are negative eigenvalues of $A$ (Fig.\,3) and continue the logarithm along that contour.

\begin{figure}[!ht]
\begin{center}
\includegraphics[height=1.3in]{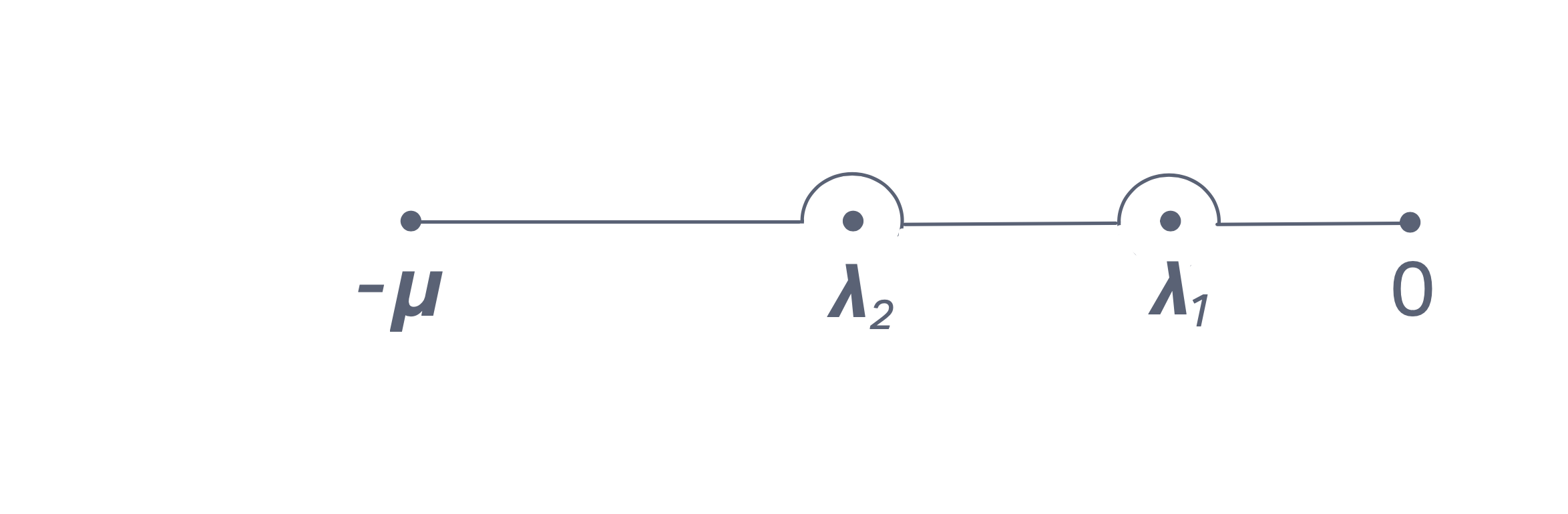}
\caption{\footnotesize The appropriate contour for the function ${\mathrm{Ln}}(E+\mu A)$ in case of negative eigenvalues.}
\end{center}
\end{figure}

\bigskip

\noindent\textbf{Proof of Lemma 1}. 

\bigskip

\noindent\textbf{1.} First of all, we consider the case $f(t)\equiv 0$ that is linear homogenuous system. Fix a $t_0 \in {\mathbb T}$. Consider the fundamental matrix $\Phi(t)$ of the time scale system 
$$x^\Delta=A(t) x$$
with initial conditions $\Phi(t_0)=E$. 
Then
\begin{equation}\label{eqord}
\begin{array}{l}
\Phi(t)=E+\int_{t_0}^t A(\tau) \Psi (\tau)\, \Delta\,\tau= \\
E+\int\limits_{[t_0,s(t)]\cap {\mathbb T}} A(\tau) \Phi (\tau)\, d\,\tau+\sum\limits_{\tau\in [t_0,t]\cap {\mathbb T}:\mu(\tau)>0} A(\tau)\Phi(\tau)\mu(\tau).
\end{array}
\end{equation}

On the other hand, for any $\tau\in {\mathbb T}$ with $\mu(\tau)>0$, we have 
\begin{equation}\label{eqphi}
\begin{array}{l}
\Phi(\sigma(\tau))=(E+\mu(\tau)A(\tau))\Phi(\tau)=
\exp(B(s(\tau))(\ln(1+\mu(\tau)))\Phi(\tau)=\\
\Phi(\tau)+\int_0^{\ln(1+\mu(\tau))} B(s(\tau)) \Phi(\tau+\theta)\,d\theta
\end{array}
\end{equation}

And now, we introduce $\Psi(t)$ as a fundamental matrix of the system 
$$\dot x=B(t)x$$
such that $\Psi(t_0)=E$.

Obviously, this matrix satisfies the equation 
\begin{equation}\label{eqpsi}
\Psi(s)=E+\int_{t_0}^s A(\tau) \Psi (\tau)\, d\,\tau
\end{equation}
Given a $t\in {\mathbb T}$, we obtain from \eqref{eqord}, \eqref{eqphi} and \eqref{eqpsi} that matrices 
$\Phi(t)$ and $\Psi(s(t))$ coincide being solutions of the same integral equation \eqref{eqord}.

Now, we apply a similar approach to a non-homogeneous system. 

Fix a $\tau$: $\mu(\tau)>0$. Let $x(t)$ be a solution of Eq,\, \eqref{e8}. By definition, 
$$x(\sigma(\tau))=(E+A(\tau)\mu(\tau))x(\tau))+f(\tau)=
\exp(\ln(1+\mu(\tau))B(s(\tau)))x(\tau)+f(\tau)\mu(\tau).$$

On the other hand, in order to obtain a solution of Eq. \eqref{e8} we must have 
$$x(\sigma(\tau))=
\exp(\ln(1+\mu(\tau))B(s(\tau)))x(\tau)+
\int_0^{\ln(1+\mu(\tau))} \exp(B(\mu(\tau)-\theta))g(s(\tau))\, d\theta$$
which implies the second equality of Eq.\, \eqref{e9}.
$\square$

Later on, considering linear systems \eqref{e7} on time scales we assume the following.

\bigskip

\noindent\textbf{Condition I.} \emph{
\begin{enumerate}
\item Matrices $A(t)$ and $(E+\mu(t) A(t))^{-1}$ are continuous and bounded on $\mathbb T$.
\item Either the time scale is syndetic or  $A(t)$ is invertible for all $t\in {\mathbb T}$ and $\limsup_{t\in {\mathbb T}} \|A^{-1}(t)\|<+\infty$.
\end{enumerate}}

\bigskip

\noindent\textbf{Remark 2}. The matrix $$A^{-1}(t_0){\mathrm{Ln}}\,(E+\mu(t_0)A(t_0))$$ can be well-defined for non-invertible matrices $A$. If $A$ has zero eigenvalues only, this matrix can be defined as
$$\sum_{k=0}^\infty (-1)^kA^k\mu^{k+1}/(k+1).$$
The sum is finite which is evident for a Jordan block. For other non-invertible matrices, the above expression may be calculated block-wise. If $A$ and $f$ are continuous on $\mathbb T$ then $B$ and $g$ are continuous on $\mathbb S$ (however, they may be discontinuous on all the axis $\mathbb R$).

\bigskip

\noindent\textbf{Lemma 2.} \emph{Let Condition I be satisfied and $f(t)$ be bounded and continuous. Then the matrix function $B(s)$ and the vector function $g(s)$ defined by Eq.\, \eqref{e9} are bounded and piece-wise continuous on the real axis.}

\bigskip

\noindent\textbf{Proof}. First of all, let us prove this statement for $B(s)$. Assume that there is a sequence $t_k$ such that 
\begin{equation}\label{equinfty}
\|B(s_k)\|\to\infty
\end{equation}
where $s_k=s(t_k)$. Let $A_k=A(t_k)$, $\mu_k=\mu(t_k)$. We may assume that the matrices $A_k$ converge to $A_0$, and $\mu_k\to \mu_0\in [0,\infty]$.

If $\mu_0\notin \{0,\infty\}$ we obtain that the matrix $E+\mu_0 A_0$ is degenerate that contradicts to our assumptions. 

Let $\mu_0=0$. We easily obtain that $\lim B(s_k)=A_0$ that contradicts to \eqref{equinfty}. 

Let $\mu_0=+\infty$. Then 
$$\begin{array}{c}
\lim
\dfrac{{\mathrm{Ln}}(E+\mu_k A_k)}{\ln(1+\mu_k)}=
\lim
\dfrac{{\mathrm{Ln}}[(\mu_k E)(E\mu_k^{-1}+A_k)]}{\ln(1+\mu(t_k))}
=\\
\lim\dfrac{{\mathrm{Ln}}[(\mu_k E)(E\mu_k^{-1}+A_k)]}{\ln(1+\mu_k)}=
\lim\dfrac{({\mathrm{Ln}}\mu_k E)+{\mathrm{ Ln}}(E\mu_k^{-1}+A_k)]}{\ln(1+\mu_k)}.
\end{array}$$
We consider two cases. If $A_0$ is non-degenerate, the above limit equals $E$. In the opposite case, we appeal to the second assumption of Condition I. $\square$

\bigskip

\noindent\textbf{Remark 3}. The last statement justifies the renormalization procedure, introduced above. Indeed, if we leave the time scale 'as is', all bounded systems become non-hyperbolic. However, non-homogeneous systems still might have bounded solutions.

\bigskip

\noindent\textbf{Definition 8.}
 We say that a system \eqref{e1} is \emph{hyperbolic} on an interval $(t_1,t_2)$ which can be unbounded, if the equivalent system 
$$ \dot x = B(t)x$$ is hyperbolic on $(s(t_1),s(t_2)$.

\bigskip

Now, let us consider the continuous family of matrices on a time scale $\mathbb T$: 
$$A(t,\nu), \qquad \nu \in {\mathbb R}^n.$$
We consider a non-homogenuous system 
\begin{equation}\label{e10}
x^\Delta = A(t,\nu)x+f(t) 
\end{equation}
with the function $f$ being uniformly bounded $|f(t)|\le 1$ and the corresponding family of homogenuous systems 
\begin{equation}\label{e11}
x^\Delta = A(t,\nu)x 
\end{equation} 
which is assumed to be uniformly hyperbolic (the constants $a$ and $\lambda$ may be selected the same) on families of segments.

\bigskip

We assume that the above Condition I is satisfied uniformly to the parameter $\nu$. However, we would like some other conditions to be fulfilled.

\bigskip

\noindent\textbf{Condition II.} There exist positive constants $a$ and $\lambda$ such that for any $\nu\in {\mathbb R}^m$ there exist 
$$-\infty=\tau_0<\tau_1<\ldots<\tau_k=+\infty$$
such that for any $\nu$ the system \eqref{e11} is hyperbolic on all segments $[\tau_{j-1},\tau_j]$, $j=1,\ldots, k$.
Here the number $k$ and all values $\tau_j$ may depend on the parameter $\nu$. 
Let $E^s_j(t,\nu)$ and $E^u_j(t,\nu)$ be the corresponding stable and unstable spaces (any, if the segment is finite and multiple choice is possible).

\bigskip

\noindent\textbf{Condition III.} $\dim  E^s_j(\tau_j,\nu)<\dim  E^s_{j+1}(\tau_j,\nu)$ for any $j$ and $\nu$. 

\bigskip

\noindent\textbf{Condition IV.} There exists a value $\alpha>0$ that does not depend on $\nu$ and such that $\angle (E^u_j(\tau_j,\nu), E^s_{j+1}(\tau_j,\nu))>\alpha$.

\bigskip
It follows from Lemma 2 that for non-syndetic time scales conditions I-IV cannot be true unless all systems are uniformly unstable hyperbolic ($E^u={\mathbb R}^n$ and the constants can be taken the same) over the real line.
\bigskip

\noindent\textbf{Theorem 2}.\emph{Let the homogeneous system \eqref{e11} satisfy conditions I-IV uniformly. Then there exists a value $K>0$ such that for any right hand side $f$ with  $|f(t)|\le 1$ there exists a solution $\varphi(t,\nu)$ of Eq.\, \eqref{e10} such that  $|\varphi(t,\nu)|\le K$ for any $\nu$. Moreover, increasing the value $K$, we may obtain a linear operator ${\mathcal L}_\nu$ such that $\varphi(\cdot,\nu)={\mathcal L_\nu}f$ and, for any fixed value $t$ the function $\varphi(t,\nu)$ continuously depends on $\nu$.}

\noindent\textbf{Proof.} Lemma 1 demonstrates that there is a continuous embedding of continuous functions on a time scale to ${\mathbb L}^\infty ({\mathbb R})$. Due to uniform hyperbolicity of system \eqref{e11}, there exists uniformly bounded systems of linear operators that map right-hand sides of ${\mathbb L}^\infty ({\mathbb R})$ to bounded solutions of corresponding linear non-homogeneous systems of ordinary differential equations.

Then, the reduction of these solutions to the time scale gives us a family of bounded solutions on the time scale.
$\square$

The proof of the above theorem follows from the reduction, introduced above and, also, from the similar result for real line obtained by V.\, A.\, Pliss \cite{P81}. 

\section{Structural stability for time scale systems}

Now we proceed to formulate a conjecture on structural stability for the time scale dynamics. 

\bigskip

Consider a system
\begin{equation}\label{e12}
x^\Delta
\bigskip
=V(t,x), \qquad t \in {\mathbb T}.
\end{equation}
We suppose that 
\begin{equation}\label{e13}
|V(t,x)|
\bigskip
<M \qquad \dfrac{\partial V}{\partial x}(t,x) <M 
\end{equation}
for some $M>0$ and that the above partial derivative is uniformly continuous. Suppose that all solutions of the considered system satisfy the forward uniqueness condition which can be provided by $rd$ -- continuity condition.
Let $\varphi(t,x_0)$ be the solution of (11) with initial conditions $x(0)=x_0$.

Consider the perturbed system 
$$x^\Delta=V(t,x)+Y(t,x); $$
with
$$|Y(t,x)|\le \delta, \qquad |D_xY(t,x)|\le \delta.$$
defined on a time scale ${\mathbb T}_1$ such that 
$$d_H({\mathbb T},{\mathbb T}_1)<\delta.$$
Here $d_H$ is the Hausdorff distance between two closed subsets of $\mathbb R$ (of course, this distance may be infinite).
Take the linearization of system \eqref{e12} in a neighborhood of the solution 
$\varphi(t,x_0)$
$$
u^\Delta=A(t,x_0)u
$$
where $$A(t,x_0)=\dfrac{\partial V}{\partial x} (t,\varphi(t,x_0)).$$

\bigskip

\noindent\textbf{Conjecture.}\emph{
Let system \eqref{e13} satisfy the above conditions I-IV and for any $j$ the conditions
$$\tau_{j+1}-\tau_j> T(a,\lambda,\alpha), \qquad j = 0,\ldots s.$$
Then for any $\varepsilon>0$ there exists a $\delta>0$ such that if condition (16) is satisfied, there exists a homeomorphism $\varphi$ of the space ${\mathbb R}^n$ such that for any $x_0 \in {\mathbb R}^n$
$$|x(t,0,x_0) - y(t,0,\varphi(x_0))| <\varepsilon.$$}

\bigskip

We hope that the results of this paper will be a keystone that allows us to prove the above conjecture.

\end{document}